# Blueprint
# for a
# Classic Proof
# of the
# 4 Colour Theorem

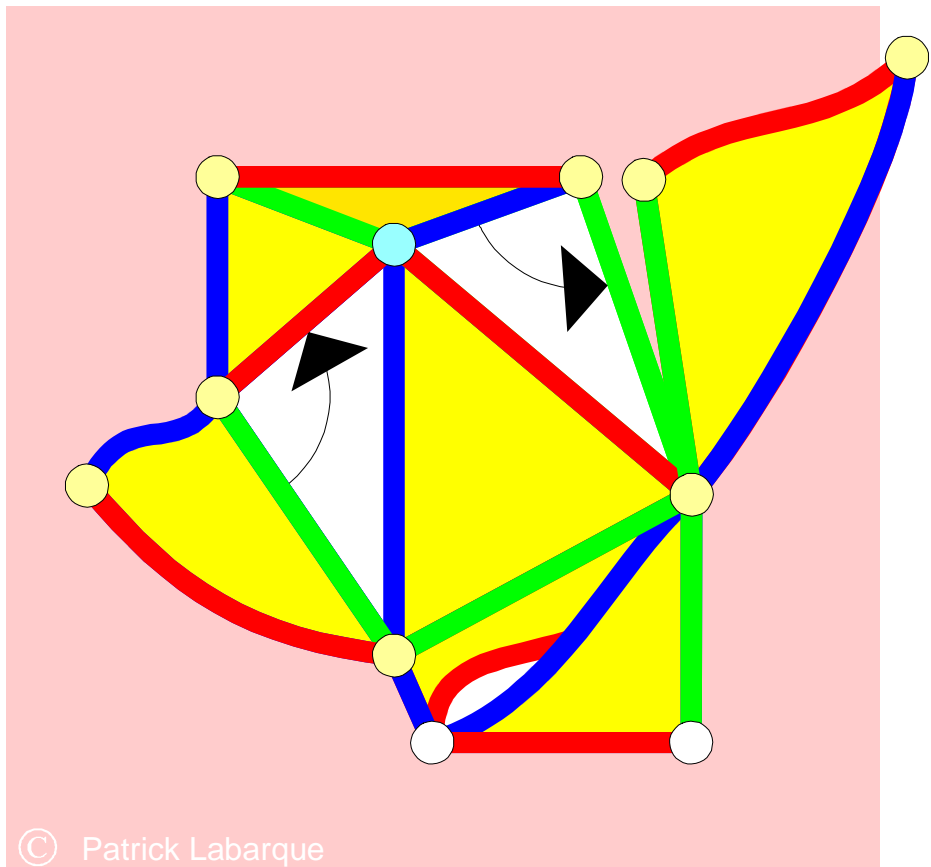

## Patrick Labarque

# Blueprint 3a for a Classic Proof of the 4 Colour Theorem

"or making stable statements in a changing world"
Patrick Labarque: patricxk.labarxque@verxsateladsl.be (delete any x)
2008 august 23

## 1   Abstract

*The proof uses the property that the vertices of a triangulated planar graph with v vertices can be four coloured if the triangles of it can be given a +1 or -1 orientation so that the sum of the triangle orientations around each vertex is a multiple of 3. Such orientation is first used separately on one of the two triangulated polygons resulting from a Hamilton circuit in a triangulated planar graph with v vertices. The graph is then reconstructed by adding the triangles of the other polygon one by one. When the graph is totally reconstructed there is always a combination for the orientations of the triangles for which their sum around each of v-2 successive vertices in the Hamilton circuit is a multiple of 3. It is then provable that the sum of the triangle orientations around the two remaining vertices must also be a multiple of 3.*

## 2   Triangulated planar graphs

The 4colour theorem was first stated as:
*Prove that 4 colours are sufficient to colour any planar map so that two countries with a common border never have the same colour.*
Instead of a planar map its dual can be used as well. One vertex is placed in each country and the vertices in adjacent countries are joined by an edge that crosses the common border. In this new graph the same question can now be asked for the adjacent vertices. It is made more difficult when a maximal number of non crossing edges are added. This is called a maximal planar graph (MPG). In such a maximal planar graph all faces are triangles (if not, it's not maximal as edges can be added in any face with more than 3 sides). Contrary any planar graph where all the faces are triangles (included the infinite face) is a MPG. Instead of MPG, we use "triangulated planar graph" hereafter as it expresses better the fact that all faces are triangles. It expresses also better the relation with triangulated planar polygons where all faces inside the polygon are also triangles.
The 4colour theorem has a lot of equivalent formulations for triangulated planar graphs.

## 3   Colour schemes for triangulated planar graphs

In the text hereafter we use "colouring" when adjacent elements must be different, and "numbering" is used when certain elements must have a numeric property (e.g. be equal to 0).

### 3.1   4-colouring of the vertices

The most common known scheme is *the 4colouring of the vertices so that no vertex is adjacent to a vertex with the same colour*. We call this a **V4c (Vertex4-colouring)**. Here we use CMYK (abbr. of **C**yan, **M**agenta, **Y**ellow and blac**K**) as the four vertex colours.

### 3.2   3-colouring of the edges

Another scheme, called a Tait colouring, is *a 3colouring of the edges so that each triangle has 3 different coloured edges* (Tait 1880). We call this an **E3c (Edge3-colouring)**. Here we use **rgb** (abbr. of **r**ed, **g**reen and **b**lue) as the three different edge colours. Instead of colours we can also use numbers. It has the advantage that one can calculate with them. By convention we say r=0, g=1 and b=2.

### 3.3 Orientation of the triangles and mod3 numbering of the vertices

One gets yet another scheme by giving an orientation (±1) to each triangle (Heawood 1898) so that their sum around each vertex is a multiple of 3. We call this a **T2# (Triangle2-numbering)**. If one makes the sum mod3 of the T2#'s adjacent to a vertex one gets a 0, 1 or 2 for that vertex. We call this a **V3# (Vertex3-numbering)**. *A good combination of T2# for all the triangles results in a V3#=0 for all the vertices*. Here 1 and 2 are used instead of +1 and -1 as it is shorter, easier to make the sum mod3 around a vertex, but especially because it's easier to make the relation with the E3c. The nice thing about a T2# is the fact that even when the sum mod3 of the T2#'s around a vertex is not 0, we get a consistent V3# for that vertex or even a partial V3# with part of the triangles around a vertex. Combining two such partial numberings different from 0 can result in a "good numbering" for the whole vertex. Because of these properties this scheme is used for the proof.

### 3.4 Mutual relations between the three colour schemes

For all "good" coloured triangulated graphs these three schemes can be linked unequivocally to each other:

**a) V4c ⇔ E3c:** The mutual relation between a V4c and an E3c is given by *fig.1 left*.

**b) E3c ⇔ T2#;V3#:** To link a T2#;V3# with an E3c we need a consistent convention for the orientation of the triangles and the partial V3# in relation with the Edge3 numbers as in *fig.1 right*. If the T2#=1 the rgb (0,1,2) order of the E3c in a triangle is clockwise, if the T2#=2 the rgb order is counter clockwise.
In a triangulation with a good E3c expressed as numbers within this convention we have:
mod3 of [the E3c of the $1^{st}$ edge + the sum of the partial V3#'s at the right side of a path] = the E3c of the last edge on that path.
E.g. with *fig.1* right with path BCDE: BC as the $1^{st}$ edge is 2 (blue), DE as the last one is then equal to mod3[E3#$_{BC}$+V3#$_{BCD}$+V3#$_{CDE}$] or mod3[2+0+1]=0, or edge DE must be red.
This applies also to two edges of the same vertex.

**c) T2# ⇔ V4c:** the mutual relation is the following:
A vertex D of a triangle ACD, adjacent to triangle ACB, has the same colour as B when both triangles have a different orientation[1], and D has the fourth colour if both orientations are the same (see *fig.1 right*).

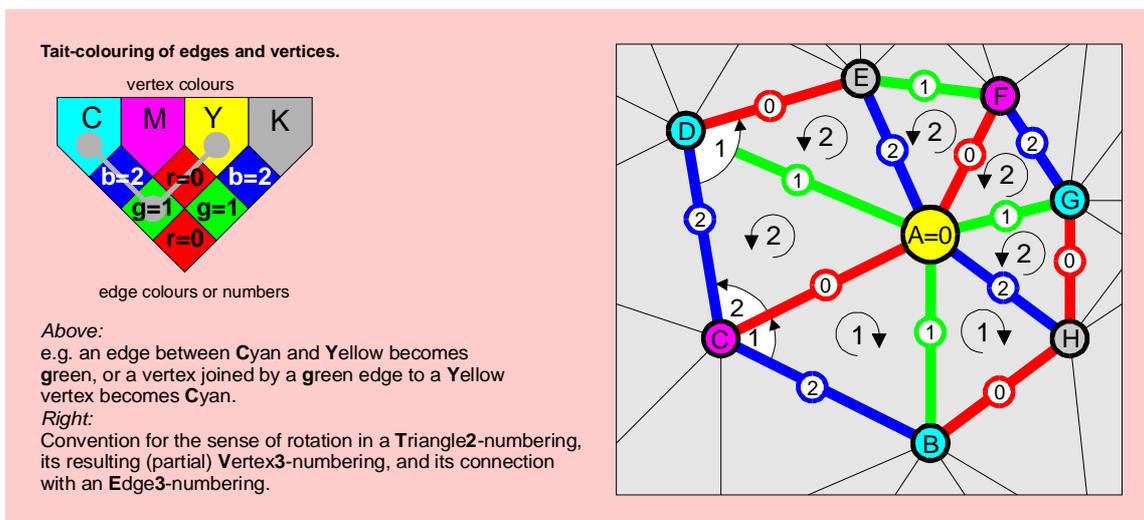

*fig.1: Convention and mutual relations between the three colour schemes.*

---
[1] Thus: if in a triangulation all adjacent triangles can have a different T2# its vertices can be 3coloured.

### 3.5 *Number of different colourings for a triangulated planar graph in the three schemes*

Hereafter we often have to speak about the combinations of numberings. We do that by putting a C before the abbreviation. e.g.: a CT2# for x triangles indicates one combination of Triangle2-numberings for the set of x triangles.

**... in a CT2#;CV3#:**

A CT2# for a triangulated planar graph has always a complementary one, by switching the 1's and 2's. It results in a complementary CV3# for all the vertices as if we kept the 0's unchanged and switched here also the 1's and the 2's. A CV3# is different from its complement, except if all V3#'s are zero.

A good CT2# results in a CV3# with all V3#=0. Its complement gives the same CV3#. Good CT2# come by pairs. Taking one CT2# from each good pair gives n morphological different good CT2# from a set of 2n different good CT2#. The easiest way to make such a set of n morphological different good CT2# is to keep one triangle unchanged.

**... in an E3c:**

Any of 2n good CT2# for a triangulated planar graph can be transposed into a CE3c using the algorithm from section 3.4.b. The resulting CE3c is then totally defined when a colour is given to the first edge. There are three possibilities, r, g or b. Thus any CE3c belongs to a set of 6 different but related CE3c. Permutation of the edge colours in one colouring produces the 5 others.

**... in a V4c:**

Any of the 6n good CE3c for a triangulated planar graph can be transposed into a CV4c using the scheme of *fig.1 left*. Giving a colour to one of the vertices defines the colours of the other vertices. We have four possibilities C, M, Y or K. Thus any CV4c belongs to a set of 24 different but related CV4c.

Thus each morphological different good "colouring" of a triangulated planar graph belongs to a set of 24 related V4c, 6 related E3c and 2 related T2#. The V3# is thus also the densest of the three colour schemes. One more reason to use it for the proof as it's easier to handle two permutations instead of six or twenty four.

## 4 Relations between Triangle2- and Vertex3-numberings

Instead of using the usual edge degree for a vertex we use here the triangle degree. As we are going to combine two triangulated polygons we can then simple add the triangle degrees for the same vertex to obtain the triangle degree of the whole vertex. This agrees better with the T2#, the V3# and the addition of two partial V3#'s of the same vertex.

### 4.1 … *for 1 triangle:*

The number of different CT2# for a triangle is $2^1$, as a triangle can only have two different T2# (1 or 2). They are complementary, as substitution of 1 by 2 gives the other T2#. The resulting partial V3# for a corner of a triangle is equal to the T2# for that triangle, it can never be zero.

### 4.2 ... *for 1 vertex:*

The number of different CT2# for t triangles around a vertex is equal to $2^t$. A vertex can have only three different V3#'s (0, 1 or 2) but with t>=2 we have more than three different CT2# for one vertex. Thus each CT2# results in V3# for a vertex, but a V3# can have more than one generating CT2#. Although we don't need them for the proof we give the following formulas:
$(2^t+2\cdot(-1)^t)/2$ = the number of different CT2# that gives a V3#=0 for a vertex with degree t.

$(2^t-1\cdot(-1)^t)/2$ = the number of different CT2# that gives a V3#=1 (or 2) for a vertex with degree t.

They give the following table[2]:

| number of triangles | 1 | 2 | 3 | 4 | 5 | 6 | 7 | ... |
|---|---|---|---|---|---|---|---|---|
| **#CT2# with V3#=0** | 0 | 2 | 2 | 6 | 10 | 22 | 42 | ... |
| **#CT2# with V3#=1** | 1 | 1 | 3 | 5 | 11 | 21 | 43 | ... |
| **#CT2# with V3#=2** | 1 | 1 | 3 | 5 | 11 | 21 | 43 | ... |
| **TOTAL #CT2#** | 2 | 4 | 8 | 16 | 32 | 64 | 128 | ... |

The V3# for a vertex can also be split into two or more partial V3#'s. This partial V3#'s can be added Mod3 in the same way as if we add the T2#'s separately. In preparation for the proof hereafter we look to the addition of two partial V3#'s adjacent to the same vertex.

0+0 gives a V3#=0   1+1 gives a V3#=2
0+1 gives a V3#=1   1+2 gives a V3#=0
0+2 gives a V3#=2   2+2 gives a V3#=1

### *4.3   ... for a triangulated planar graph:*

For a triangulated graph with t triangles we have $2^t$ different CT2#. As we have 2*(v-2) triangles in a triangulated graph we have then $2^{2(v-2)}$ different CT2#.

As the sum of the V3#'s of a triangulated graph must be a multiple of 3 (each corner of a triangle is counted three times) the number of different CV3# for a triangulated graph can maximal be equal to $3^{v-1}$. With v>=6 the last one is less than the number of different CT2# and some CT2# must result in the same CV3#. Thus each CT2# results in one CV3#, but a CV3# can have more than one generating CT2#.

But remember a V3# cannot exist on its own; it is always the result of a T2#. There are CV3# that cannot exist in any triangulated planar graph. As a corollary of theorem 4 hereafter, a CV3# with two adjacent V3#'s of 1 and 2 and the others equal to zero, cannot be the result of a CT2# (even if their sum is a multiple of 3). Depending on the constellation of the graph there are also other CV3# that cannot exist. As a corollary of theorem 3 hereafter the number of different CV3# for a triangulated planar graph[3] is minimum $3^{v-2}$.

For a vertex, one can calculate the number of different CT2# that gives a V3#=0, 1 or 2. Up to now there is no easy way to calculate the number of different CT2# with V3#=0 for all the vertices. If there was one, it should be a proof of the 4colour theorem.

But to prove the 4colour theorem it's sufficient to prove that between all the CT2# for a triangulated planar graph there is at least one pair with a V3#=0 for all the vertices.[4].

## 5   Triangulated planar graphs and Hamilton circuits

Any triangulated graph without separating triangles[5] has a Hamilton circuit[6] (Whitney 1931). Moreover any of the edges can be part of such a circuit (Thomassen 1983). This circuit divides the graph into two triangulated polygons; an outer one and an inner one, with all the vertices on the common perimeter (see *fig.2*).

---

[2] See SLOANE's "On-Line Encyclopedia of Integers" sequences A078008 and A001045
[3] Graphs exist with this minimum, e.g. the octahedron. Replacing successively triangles in such octahedral graph with another octahedron produces graphs with this property (all vertices have even degree). CONJECTURE 1: All graphs with vertices of even degree have exactly $3^{v-2}$ different CV3#.
[4] Some special cases are easy to prove:
If the degrees of all the vertices in a triangulated graph are a multiple of 3 (then v is even), giving the same orientation to all the triangles gives a V3#=0 for all the vertices and they need 4 colours in a V4c.
If all degrees are even then we can give a different orientation to any pair of adjacent triangles. This results also in a V3#=0 for all the vertices and they can be 3coloured.
[5] A separating triangle is a circuit with three edges, with vertices at both sides of the circuit.
[6] A Hamilton circuit visits each vertex once, and returns to the start. This is a sufficient condition, not a necessary one, as graphs exist with separating triangles having Hamilton circuits (e.g; the graph of fig.2).

This property is used for the proof. A Hamiltonian graph is split into two triangulated polygons, an inner one and an outer one. The graph is then reconstructed from the outer polygon by adding step by step a triangle from the inner polygon.

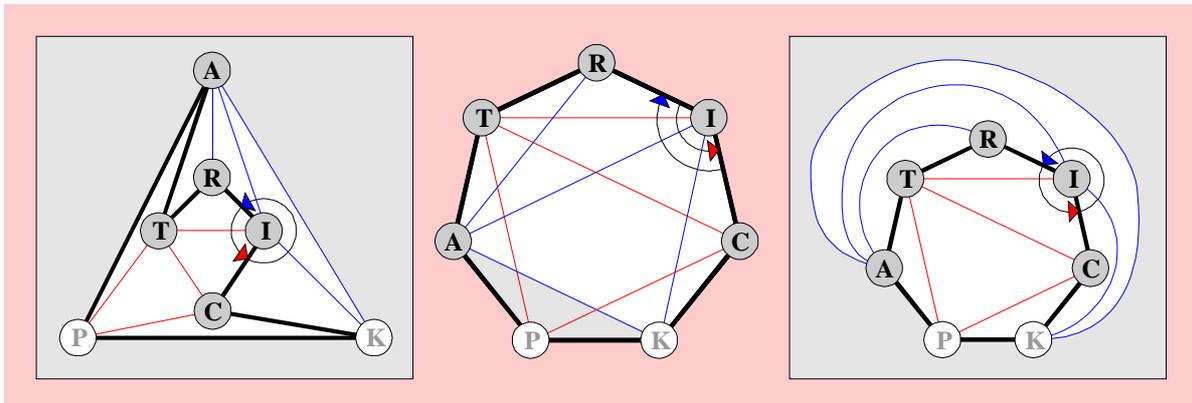

*fig.2: Left: a triangulated graph with a Hamilton circuit in it (bold). The circuit divides the graph into two triangulated polygons, an inner one (red diagonals) and an outer one (blue diagonals) with the "infinite" triangle APK (grey). Right: the same graph with the Hamilton circuit as a regular polygon. The diagonals of the outer polygon need to be curved. It allows representing triangulated graphs with double edges. Centre: the two polygons on top of each other. It can be seen as a transparent globe with the Hamilton circuit as the equator. The triangle APK (also in grey) is a normal triangle. The outer face here is not a part of the graph. Note that the outer polygon is mirrored and that the convention for the orientations must be inverted for that polygon (see vertex I).*

The use of triangulated polygons is essential for the proof hereafter, so we look first to some of their properties.

## 6 Triangulated planar polygons (without inner vertices)

### *6.1 Some form properties of triangulated planar polygons*

**a)** If one makes no difference between symmetries a polygon with 6 vertices can be triangulated in 3 ways (the red, green and blue one in *fig.2*). Making a difference between mirror images but not between rotations gives here one more (the magenta one). Making differences between rotations gives here 14 different types on a base[7], the 4 types plus all their rotations resulting in a *different* polygon.

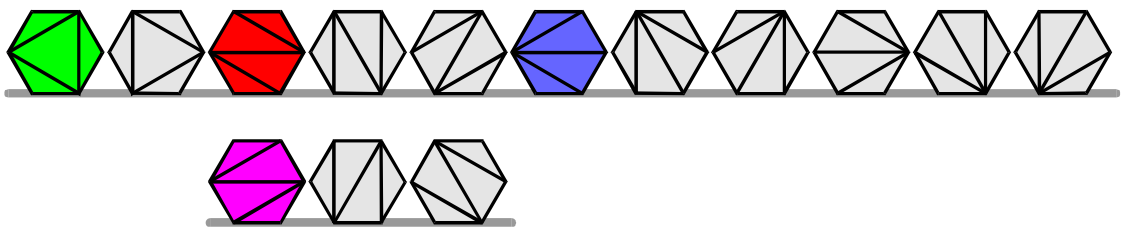

*fig.3: The fourteen different types of triangulated hexagons on a base. Any hexagon has at least two vertices with a triangle degree=1.*

---

[7] Use following algorithm to construct a genealogical tree of triangulations. Suppose we have the different triangulations with v vertices, where the "last added ear" is marked. All edges are free except the base and the edges in clockwise direction between the base and the "last added ear". Add now in each triangulation one ear to and for each free edge. Mark that ear as "last added ear" in each new triangulation with v+1 vertices. Start with the base (v=2), add one triangle and mark it as "last added ear"; continue, using the algorithm.

The number of different types of triangulated polygons on a base is given by the Catalan number. It's a function of the number of vertices. Starting with v=2 (a degenerated polygon) the first Catalan numbers are: 1, 1, 2, 5, 14, 42, 132...

**b)** Any triangulated planar graph with a Hamilton circuit gives two coupled triangulated polygons. The inverse is true as well. The set of combinations by two of the different triangulated polygons on a base contains all the triangulated planar graphs on a base with a Hamilton circuit. If they have n common diagonals the resulting graph will have n double edges. But even combining two idem triangulated polygons can still give a triangulated graph with as much double edges as we have diagonals. Such a graph is also suitable for the proof[8]. In that case *fig.2 right* is the best representation for it.

**c)** In a triangulated polygon we have three different types of triangles:
- triangles with two edges on the perimeter, called "ears". The triangle degree of the vertex at the ear tip is 1
- triangles with one edge on the perimeter (we can call them "cheeks")
- triangles with zero edges on the perimeter (we can call them "noses")

A triangulated polygon has always at least two ears (Gary H. Meister 1975) and two ear tips cannot be adjacent. Thus a triangulated polygon on a base has always at least one ear tip that is not a base vertex.

**d)** A triangulated polygon on a base has a hierarchical structure for its triangles. Starting with the "ears" each triangle has one "parent" except the last one on the base. As a corollary of this any vertex that is not a base vertex can be associated with the triangle that is most near to the base. We call it the associated vertex-triangle, and we call its associated vertex the triangle-vertex[9]. Any triangle has exactly 1 triangle-vertex and any vertex has exactly 1 vertex-triangle except the two base vertices. This is self evident as there are only v-2 triangles in a triangulated planar polygon. Using this association we give hereafter the same letter (lower case) to a triangle as its triangle-vertex (capitals). We add an ' or " to differentiate between the triangles of the inner and the outer polygon. There are no triangles with the same letter as the two base-vertices.

Be aware that choosing another base for the same triangulation gives another association and another notation for the triangles.

### 6.2 Colouring triangulated planar polygons

It is very easy to apply a good V4c or E3c to a triangulated polygon without vertices in it. Simply start to colour the elements (vertices or edges) of one triangle and then colour successively the elements of any adjacent triangle. Each of these colourings results in a CV3# for the polygon and the V3#'s don't need to be zero. They can also be seen as combinations of partial V3#'s for a triangulated planar graph where the triangles of the other polygon are not yet involved. Hereafter we use X' (or X") to indicate the partial V3# of vertex X of the inner (or outer) polygon.

THEOREM 1: ***A triangulated planar polygon with v vertices has $2^{v-2}$ different CV3#, one for each CT2#.***

Given a triangulated planar polygon and one of the CV3# resulting from one of the CT2#. It has at least two ear tips and the T2# of each ear is equal to the V3# (1 or 2) of its ear tip. One of the ears can be cut off and its T2# can be subtracted from the V3# of its two neighbours, to adjust the CV3# of the remaining polygon. Again this polygon has at least one ear tip and the same process can be repeated. In each step the T2# of the cut off ear was defined by the adjusted V3#. We end with one triangle and in this step the three remaining vertices have an adjusted V3# equal to the T2# of that triangle. Thus any CV3# in a triangulated polygon has ONE generating CT2#. As we have $2^{v-2}$ different CT2# we have then also $2^{v-2}$ different CV3#, one for each CT2#.

---

[8] All vertices in such a graph have an even degree. Its vertices can be 3coloured.
[9] We use "triangle-vertices" also to indicate the vertices of a triangulated planar polygon except the two base vertices.

**Corollary 1:**
A triangulated planar polygon has always an ear tip that is not a base vertex. The same process can be repeated with a polygon on a base, using in each step an ear tip that is not a base vertex. In the last step all T2# are defined thus the V3# of the two base vertices are also defined. Their V3#'s are not needed, and any edge can be chosen as a base, thus: **the CT2# and the CV3# for a triangulated polygon with v vertices is defined if the CV3# of only v-2 successive vertices is known.** This suggests that the V3#'s of two adjacent vertices may be disregarded. For the proof hereafter this observation is used by disregarding until the end the V3#'s of the two base vertices.

**Corollary 2:**
During the process in corollary 1, we see that the cut off ear and the cut off ear tip in each step are in fact the vertex-triangle and its triangle-vertex. Making equations using the same hierarchical order by which the ears are cut off we see that the partial V3# of each triangle-vertex can be expressed as a function of its vertex-triangle and the V3#'s of some of the preceding triangle-vertices (see *fig.4*). Thus for any combination of the V3#'s of the preceding vertices any triangle-vertex has two different V3#'s. It's clear that it is also the case for the other triangle-vertices that are not adjacent to its vertex-triangle. Thus: **any triangle-vertex has two different V3#'s for any CV3# of the vertices that are not adjacent to its vertex-triangle.**

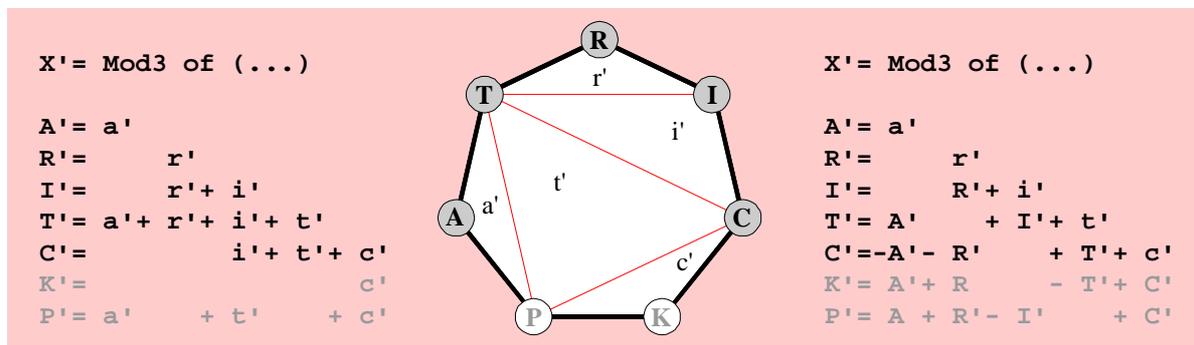

*fig.4: Centre: the inner polygon from fig.2 with the vertex-triangles marked with the same letter as their triangle-vertices. Left: equations in hierarchical order for the partial V3# of the vertices in function of the adjacent triangles. Right: equations in hierarchical order for the partial V3# of each vertex as a function of its vertex-triangle and the preceding vertices. The V3#'s of the base vertices P and K can be expressed as a function of the preceding vertices only. The equations must be seen as modular equations. This can be done for any triangulated planar polygon.*

In fact each triangle-vertex has also two different V3#'s in combination with the two vertices adjacent to its vertex-triangle but the V3#'s of this two vertices are changing together with a change in the T2# of that vertex-triangle. It's confusing when things are changing behind our back. To build further on this property we have to find some stable elements in it.

### 6.3  *Oriented pairs in the set of CV3# of a triangulated planar polygon*
**a) The X1-X2 code for triangulated polygons.**

Any V3# of a triangle-vertex can be associated with the T2# of its vertex-triangle and we can associate a 1 or 2 with each triangle-vertex. X'1 indicates then the V3# of vertex X' with a T2#=1 for its vertex-triangle x'. We can mark now a CV3# of v-2 successive vertices as a combination of Xone's and Xtwo's for each of the triangle-vertices X.

We look to the following pair from *fig.5*:

**B'1**, C'1, D'2 indicates ONE of the combinations in the polygon and
**B'2**, C'1, D'2 indicates another one where only triangle **b'** has changed.

It does NOT say that the V3# of X'1 or X'2 is equal to 1 or 2. That's only the case for an ear tip. It does NOT say that the V3#'s of the other triangle-vertices C and D are the same in both combinations (see D' in *fig.5*). That's only the case for a triangle-vertex with its vertex-triangle on the base. The ONLY thing it says is that the V3# of B' in the second combination is one more mod3 than in the first one. B'1 and B'2 form an **oriented pair** of partial V3# in any two combinations with the same X1-X2 code for the other vertices. The same is true for the oriented pairs in the other triangle-vertices. This property, reflected in the X1-X2 code, is also valuable for any other polygon independent of its constellation.

**b) The invariant structure of oriented pairs in a triangulated polygon**

The structure of this oriented pairs in the CV3# is illustrated in *fig.5* for a polygon with five vertices. It's a minimal example that contains the 3 cases when we change the T2# of a triangle (d', b' and c' with a change in respectively none, one or two *other triangle-vertices*).

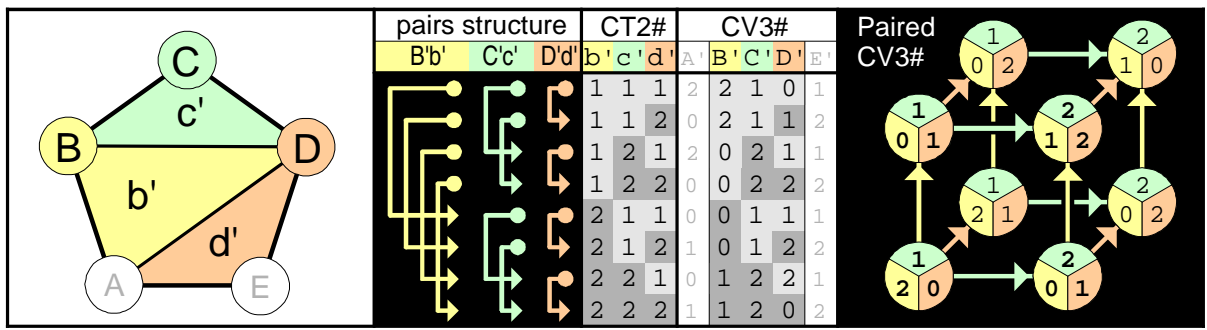

*fig.5: Oriented pairs for the polygon on the left. Black column: The arrows indicate the orientation of the pairs. The structure of the oriented pairs has 3 (=v-2) dimensions for this polygon, one for each triangle or triangle-vertex, and each dimension has four (=$2^{v-3}$) oriented pairs. There are 8 (=$2^{v-2}$) different combinations of tails and arrows, one for each row. Grey columns: The oriented pairs in the CT2# and the CV3# have the same structure. Right: representation of the oriented pairs in the CV3#'s as a three-dimensional cube (see the appendix for a visualisation with more triangle-vertices). Such schemes can be made with another choice of the base and/or for another triangulated polygon with any number of vertices.*

With *fig.5* and its caption as argumentation we see that the properties of the oriented pairs are generated by the binary structure of the CT2#'s. Because of the association between vertex-triangles and their triangle-vertices these properties are reflected in the corresponding CV3#:

1) Each vertex-triangle has $2^{v-3}$ oriented pairs of CT2#. The T2# of the tail and the arrow has a difference of 1. Each associated triangle-vertex has the same $2^{v-3}$ oriented pairs of CV3#. The V3# of the tail and the arrow of the associated triangle-vertex of a pair differ by 1mod3.
2) We define the combinations of oriented pairs as combinations of their tails and arrows. The $2^{v-2}$ different CT2# can then be expressed as combinations of oriented pairs of T2#'s for the vertex-triangles. The resulting CV3# of the triangle-vertices can also be expressed as combinations of the same oriented pairs of V3#'s for the triangle-vertices.
3) The structure of the oriented pairs is very stable.

a) Changing the T2# of a vertex-triangle changes the V3#'s of its triangle-vertex. It changes also the V3#'s of at most two other triangle-vertices, but it changes neither the structure nor the orientation of the oriented pairs in these other vertices.

b) Even if we add mod3 a number to *all* the V3#'s of one or more triangle-vertices in the set of CV3# the structure and properties of the oriented pairs are kept the same. ***The structure and properties of oriented pairs are invariant under addition*** (and subtraction).

These properties of oriented pairs are in fact more than we need. We don't need the orientation. For the proof hereafter it's sufficient to know that the tail and arrow of a pair is *different*. It is a weakening, but it has the advantage that we don't have to worry about which part is tail or arrow. We can formulate the following generalisation:

*THEOREM 2: Any triangle-vertex of a triangulated polygon has a pair of different V3#'s for each combination of the pairs in the other triangle-vertices. This difference property does not change if we add the same number mod3 to the V3#'s of one or more vertices in the set of the CV3#.*

This theorem is the key statement for the proof[10].

## 7   Colouring triangulated planar polygons with inner vertices

### 7.1   *V3# of triangulated planar polygons with inner vertices*

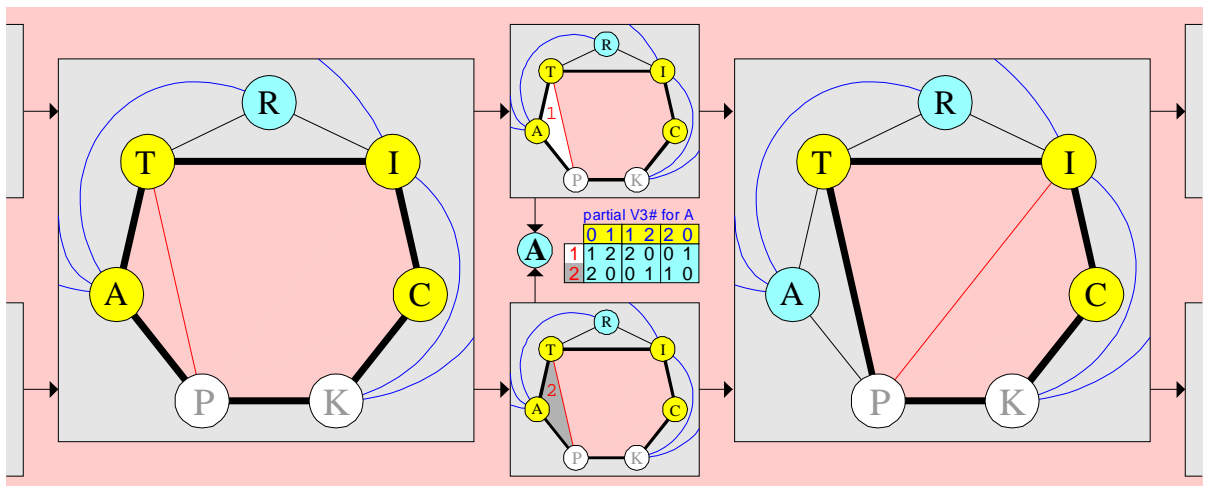

*fig.6: Hierarchical step by step addition of triangles to an outer polygon. Left: the inner cyan vertex R has three different V3#'s for any combination of the pairs of the yellow vertices on the perimeter. Any yellow vertex has two different V3#'s for any value of the trio of the cyan vertex. Centre: the triangle added to A can have a T2#=1 or 2. The four combinations with each of the three possible pairs of V3#'s for A are shown in the table in the middle; the three different V3#'s for A come up in any of these pairs. Right: vertex A becomes now a cyan inner vertex of the changed outer polygon and has now also three different V3#'s for any combination of the trio's and pairs of the other cyan and yellow vertices. This process can be repeated until no yellow vertices are left over.*

Given: an outer polygon on a base without vertices in it, and all the $2^{v-2}$ different CV3# for its v-2 triangle-vertices (again the two base vertices are disregarded). We add a triangle in the inner polygon to a triangle-vertex of this outer polygon, by joining its two neighbours. That triangle-vertex has $2^{v-3}$ pairs of different V3#'s. The triangle can have a T2#=1 or 2. We get then for each pair four V3#'s. Two of them are the same but three are always different. Each pair for that vertex is now replaced by a trio of three different V3#'s (and one double). This vertex becomes an inner vertex of the outer polygon and we adjust the perimeter of the outer polygon. The V3#s of the two other vertices adjacent to the added triangle have changed, but because of theorem 2 this does not change the difference property in the pairs of the remaining triangle-vertices on the perimeter of the outer polygon. We have now at least $3^1*2^{v-2-1}$ different CV3# in the set of $2^{v-2+1}$ generating CT2# in the new constellation.

---

[10] It is also the most vulnerable theorem; we need to "read" the figures as proof; the argumentation is not optimal, it needs a more rigid mathematical proof.

We repeat the same process (see *fig.6*) and after each step one more triangle-vertex becomes an inner vertex. Again each pair of the newly enclosed vertex is replaced by a trio and after the $n^{th}$ step the n inner vertices have $3^n$ different CV3# and for each of them any of the v-2-n remaining other vertices have a pair of different V3#'s. We have then at least $3^n*2^{v-2-n}$ different CV3#'s in the set of $2^{v-2+n}$ generating CT2#. We can thus state the following

*THEOREM 3*: **A triangulated planar polygon with $v_p$ vertices on the perimeter and $v_i$ vertices in it has all the different combinations of CV3# for the $v_i$ vertices[11]. For each of these combinations any of $v_p$-2 successive vertices has a pair of different V3#'s. The total number of different CV3# for these $v_i+v_p$-2 vertices is at least $3^{vi}*2^{vp-2}$.**

We may say "successive" as any edge of the polygon could have been chosen as the base edge. Note that we have then $2^{vp+2vi-2}=4^{vi}*2^{vp-2}$ CT2# and not all the CV3# resulting from this CT2# are different (as each trio has a double).

## 7.2   E3c of triangulated planar polygons with inner vertices

*THEOREM 4:* **In a good Edge3-coloured triangulated polygon with (or without) vertices in it, the number of edges on the perimeter with the same colour has the same parity as the total number of edges on the perimeter.**

In a good E3c of a triangulated polygon with or without vertices in it, each inner edge with colour x is adjacent to 2 triangles and each outer edge with colour x is adjacent to 1 inner triangle (and the infinite face), thus:

(1)     $t_i=2*e_i(x) + e_p(x)$

with $t_i$ the number of inner triangles, $e_i(x)$ the number of inner edges and $e_p(x)$ the number of edges on the perimeter, both with colour x. *The number of edges on the perimeter with the same colour has thus the same parity as the number of inner triangles.*

For the same polygon with Euler we have:

(2)     v-e+f=2,

with v the number of vertices, e the number of edges and f the total number of faces (included the infinite face). For the total number of faces and vertices we have:

(3)     $f=t_i+1$       (4)     $v=v_i+v_p$

with $v_i$ the number of inner vertices and $v_p$ the number of vertices on the perimeter.

If we add the number of edges of each face, every edge is counted twice thus we have $2e=3t_i+e_p$. And as $e_p=v_p$ we have:

(5)     $2e=3t_i+v_p$

Substitution of v, f and e in (2) yields:

(6)     $t_i=v_p+2(v_i-1)$

We may conclude that "*the number of inner triangles has the same parity as the total number of edges on the perimeter*".

Thus we also may conclude that "*the number of edges on the perimeter with the same colour has the same parity as the total number of edges on the perimeter*".

---

[11] We can also conclude that any inner vertex can have a V3#=0. But the theorem is more general now.

# 8  Proof

As mentioned in section 4 we first split a triangulated planar graph into two triangulated polygons on a common base, an inner one and an outer one. We start now with the outer polygon and add one by one an ear from the inner one to a non base vertex of the outer one. This can always be done as we have always an ear tip that is not a base vertex. We continue until all non base vertices have been used. We become now a "degenerated" triangulated polygon with a double edge on the two base vertices and with v-2 inner vertices[12]. Following theorem 3 the inner vertices have all the $3^{v-2}$ different combinations of V3#'s for the $2^{2(v-2)}$ CT2#. Then there is also at least one complementary pair of combinations with all inner vertices having a V3#=0.

But we still have the two vertices P and K. As all inner edges have a V3#=0 we can transpose the V3# of the triangulated polygon to an E3c, and we may colour one of the edges, say an edge on the perimeter, in any colour, say red. Each triangle inside the outer polygon has now three different edge colours. All edges have a colour, also the other edge on the perimeter, but we don't know which colour. If it's not red we cannot reconstruct the whole graph with a good E3c.

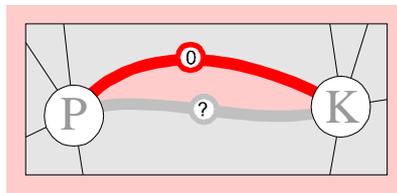

*fig.7: Part of the triangulated polygon on two vertices. The colour of the grey edge is unknown.*

But we because of theorem 4 the number of edges on the perimeter with the same colour must have the same parity as the total number of edges on the perimeter. Thus the other edge on the perimeter must also be red. As they have the same colour we can join now these two edges to one base edge and finish the reconstruction with a good E3c. Vertices P and K must have now also a V3#=0.

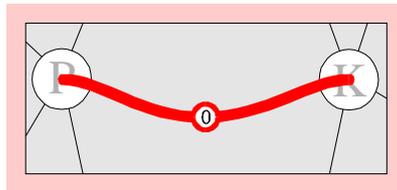

This finishes the proof. It's valuable for all triangulated planar graphs[13].

---

[12] If you don't like degenerated polygons we can also stop on the last but one triangle-vertex. We can continue the same way if the remaining polygon is a triangle (but not if it's a quadrangle).
[13] If the graph has no Hamilton circuit, then it must have separating triangles. In that case we hierarchically take out the separating triangles and their inner vertices and edges and colour the graphs separately. Then we put them back in reversed order, accept the colours of the parent and recolor the child(s) with a permutation.

# 9 Appendix: X1-X2 code and hyper cubes in triangulated planar polygons.

Let's have a triangulated polygon with three vertices on a base. $V_1$ is the ear tip and has one oriented pair of V3#'s (**V'$_1$1** and **V'$_1$2**). We add somewhere (say left of $V_1$) a new triangle $v_2$ with its ear tip $V_2$. For each oriented pair in $V_1$ we have two T2#'s for its triangle-vertex. Expressed as X1-X2 codes we have now the following CV3# for the triangle-vertices.

**V$_2$'1V'$_1$1**   **V$_2$'1V'$_1$2**   and
**V$_2$'2V'$_1$1**   **V$_2$'2V'$_1$2**

But we have an oriented pair for each vertex with the same X1-X2 code for the other vertices. We added in fact also another oriented pair to $V_1$. We have thus four pairs, a red and magenta one for $V_1$ and a blue and cyan one for $V_2$ as illustrated below.

**V$_2$'1V'$_1$1**   **V$_2$'1V'$_1$2**
**V$_2$'2V'$_1$1**   **V$_2$'2V'$_1$2**

The polygon has now two pairs for each triangle-vertex and four combinations of oriented pairs. We can continue adding vertices to the preceding combinations and each time the number of oriented pairs and combinations is doubled. After each step we have $2^{v-3}$ oriented pairs for each vertex and $2^{v-2}$ combinations of these pairs. Another way to show such expansion of a triangulated polygon is with a (v-2)-dimensional hypercube as in *fig.8* where each line represents one oriented pair. We *see* here the oriented pairs for each vertex, the combinations of these pairs, the binary structure of the X1-X2 codes and also the complement structure in the pairs and their combinations. As no specific polygon is used, it illustrates also that these properties are valuable for all triangulated planar polygons. To add a vertex we "simply" add a translated copy of the existing hypercube one unity higher and $2^{v-3}$ unities to the right; then we join all the existing dots with the corresponding new copied dots. To show the relation with *fig.5* the same colours and vertex names are used for the first three vertices.

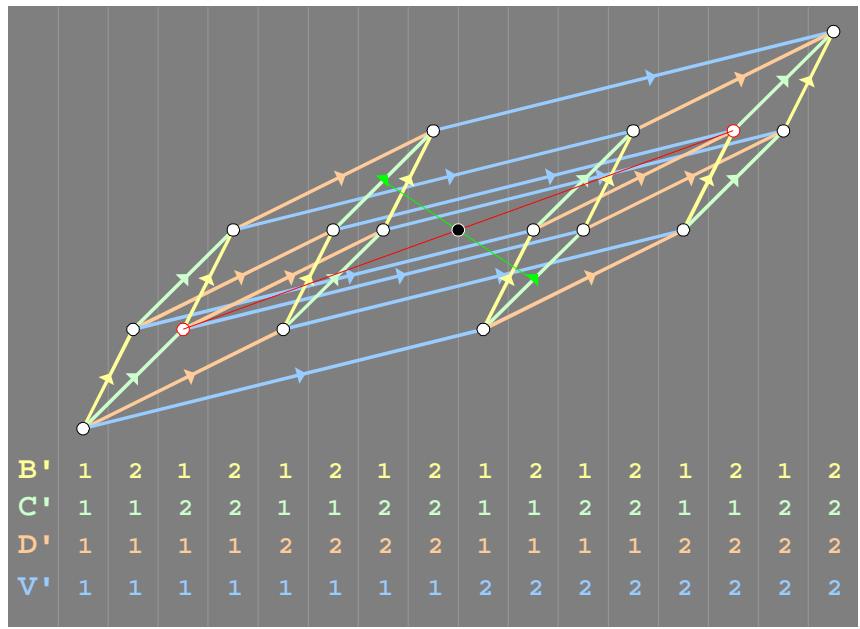

*fig.8: A 4-dimensional hypercube for a polygon with 6 vertices. The lines in the same colour represent the $2^{v-3}$ oriented pairs of V3#'s for one vertex. The arrows indicate the orientation. The white dots indicate a combination of tails and arrows (differences) of these pairs. The X1-X2 code is indicated for each combination. Each pair (e.g. the green line) and each combination (e.g. the red line) has a complementary one, symmetric with respect to the central black dot.*

# 10 Epilogue

1) I want to thank Hendrik Van Maldeghem from the RUG (Rijks Universiteit Gent, Belgium) for his valuable remarks on the first and the second proof. He e-mailed me, short but to the point, why the first and the second proof were not OK. Without his remarks this improved proof should not exist.

2) I also want to thank Christopher Carl Heckman from Arizona State University USA for his remarks on the 3$^{rd}$ version. They helped me to find and correct some stupidities.

3) I tried to make the proof as accessible as possible for everybody. Nevertheless there is still too much encoding (V4c, E3c, T2#, V3#, CT2#, CV3# etc...) in it. Sorry (see comments in the bibliography).

4) I don't master very well the mathematical language (idem for English). The proof is hybrid. It is neither a rigid mathematical proof nor a popular version of the rigid mathematical one that has to come. It can be improved and simplified in both ways. If someone sees the possibility please do it, and don't worry about authorship from my part. Anything in this paper, except the front page illustration may be used (and/or misused) by anybody.

4) The word "classic" is used in the title. I should have preferred to use the word "elegant" or "simple" but I can't because the existence theorem of Hamilton circuits by Whitney "...consists of a lengthy enumeration of cases..." ([1] page 111).

5) During research (amusement) I switched often between logic and making inventories. For the last one Excel is a marvellous tool. This program had an essential influence to exclude bad thinking and/or to discover hidden structures and relations (e.g.: theorem 3, applied to a triangulated graph without the base vertices, was discovered with limited inventories before the proof was found).

6) The first "Blueprint proof" was winter work, the result of being retired for about 4 months, and having the possibility to do what I want, without interruption. Some improvements have been made 15 months later. But this first proof and the second one were not complete. I hope this time it is OK. I still enjoyed the time spent to it and once again Lieve didn't enjoy it as much as I did. Once again, thanks Lieve and once again I promise you not to find too quickly another problem.

# 11 Bibliography

[1] Thomas L.Saaty and Paul C.Kainen "The Four-Color Problem" Dover publications 1986, ISBN: 0-486-65092-8

Nearly all stuff that is used here is in that book. It's a very complete work and low budget. Nevertheless as a layman, I needed the whole week to decipher the notations of what I had read during the weekend. And I believe it's even not the worst. What I mean: I love popularisation and especially go betweens between amateurism and academic research.

[2] Shalom Eliahou and Cédric Lecouvey "Signed permutations and the four color theorem" 2006 on the internet site: arXiv.math/0606726v1.

Similar constructions with triangulated polygons are used in this paper in a more rigid mathematical way.

[3] Additional information about triangulated planar polygons and Hamilton circuits comes from different internet sites. Unfortunately for most of the original sources you need a password[14]. So I could not check these.

Because of this lack of information I only mention the author of a theorem and the year when it was proved.

But lucky we, there is also the open internet.
There are a lot of informative websites about the subject and mathematics in general. E.g. the website of:
- Wikipedia
- Wolfram Mathworld
- The monthly Feature Column of the AMS (e.g. have a look at the marvellous animations in the column of Nov. 2006 by Jos Leys and Etienne Ghys on strange attractors, knots and lattices).
- And many others, as websites on Venn-diagrams, teasers for young would-be mathematicians, whole academic courses as "Diestel graph theory", Dr Math where one can ask questions and gets an answer[15]..., many many others that I don't see the possibility to mention them all.

---

[14] I don't blame the editors for this as there are a lot of very courageous editors.
[15] I once asked a question about the existence of Pythagorean rectangular blocks where ALL distances between the corners are integers. I got an answer back, short but to the point, so to the point that I gave up doing research on it. It's still an open problem. Some call it a perfect Euler brick or perfect cuboid.